\documentclass[11pt]{article}
\usepackage[margin=1in]{geometry}                
\usepackage{graphicx}
\usepackage{amssymb}
\usepackage{epstopdf}
\usepackage{amsthm}
\usepackage{amsmath}
\usepackage{tikz} 
\usepackage{epigraph}
\usepackage{txfonts}
\usepackage{endnotes}
\usepackage{bbold}

\usepackage{perpage}
\MakePerPage{footnote}
\tikzset{node distance=2cm, auto}
\usepackage{calc}
\newlength{\figlabelwidth} 
\newlength{\imgwidth} 

\title{Categories within the Foundation of Mathematics}
\author{Benjamin Horowitz}
\date{May 3, 2013}                                           

\begin{document}
\maketitle

\epigraph{ÒThe language of categories is affectionately known as ``abstract nonsense," so named by Norman Steenrod. This term is essentially accurate and not necessarily derogatory: categories refer to ``nonsense" in the sense that they are all about the ``structure," and not about the ``meaning," of what they represent.Ó } {Paolo Aluffi (1960-Current)\\ Algebra: Chapter 0}

\section{Introduction}

Category theory has proven to be an extremely useful tool in understanding current deep questions in mathematics. To a mathematician, category theory allows one to organize mathematical experience and ``think bigger thoughts."\endnote{Abramsky, Samson, and Nikos Tzevelekos. ``Introduction to categories and categorical logic." \emph{New Structures for Physics.} Springer Berlin Heidelberg, 2011. 3-94.} To a philosopher, category theory offers a fresh perspective on structuralist foundations to mathematics and a well developed alternative to the focus on set theoretic constructions, like the traditional type theory of Russel and Whitehead.

The recent trend in mathematics is towards a framework of abstract mathematical objects, rather than the more concrete approach of explicitly defining elements which objects were thought to consist of. A natural question to raise is whether this sort of abstract approach advocated for by Lawvere, among others, is foundational in the sense that it provides a unified, universal, system of first order axioms in which we can define the usual mathematical objects and prove their usual properties. In this way, we view the ``foundation" as something without any necessary justification or starting point. Some of the main arguments for categories as such a structure are laid out by MacLane as he argues that the set-theoretic constructions are \emph{inappropriate} for current mathematics as practices, and that they are \emph{inadequate} to properly encompass category theory itself and therefore cannot properly encompass all of mathematics, while category theory can be used to describe set theory and all the natural consequences of a given primitive system.\footnote{We will generally refer to a ZFC system unless otherwise specified.}

Others view this issue quite differently as a matter between understanding \emph{organization} verses \emph{foundations} of a mathematical theory. Kreisel argues that the category theory is, at its heart, concerned with efficiency and is therefore an organizational system, while set theory maintains its true role as arbiter of validity and is therefore a foundational system.

In this paper, we begin by describing the two basic elements of category theory; the category, a collection of objects and arrows, and the functor, a mapping between categories. From there we develop the  theory of the closed Cartesian category and the adjoint relation. This construction is at the core of Lawvere's categorical logic and is also the source of much criticism by more ``traditional" logicians, like Kreisel and Halpern. We will explore these criticisms as well as the relationship between these theories and work towards an argument for the importance of category theory as a foundational tool in mathematics. 

\section{Construction of Category Theory}

This section is meant to provide the un-initiated an overview of the key results of category theory in relation to mathematical logic and provide the terminology for much of the rest of the paper. Specific definitions are not overly important, but it is essential to get a flavor for the structures that will play a role in categorical logic.

\subsection{Categories}

A category, at its core, is a collection of objects and directional relationships between them. A collection itself is a mathematical object that contains other mathematical objects. Its relationship to a ``set" is not specified \emph{a priori}, although in most applications it will effectively function as such. 

\newtheorem{mydef}{Definition}
\newtheorem{myprop}{Proposition}
\newtheorem{myproof}{Proof}
\newtheorem{myth}{Theorem}
\newtheorem{myex}{Example}
\begin{mydef}
A \textbf{category} $\mathcal{C}$ consists of the following:
 \begin{itemize}
 \item A collection $O(\mathcal{C})$ of \textbf{objects} which are denoted by $A$, $B$, $C$, etc.
 \item A collection $A(\mathcal{C})$ of \textbf{arrows} which are denoted by $f$, $g$, $h$, etc.
 \item Two unique mappings \emph{dom, cod} $: A(\mathcal{C}) \rightarrow O(\mathcal{C})$ which assign to each arrow of $A(\mathcal{C})$ its \textbf{domain}, \emph{dom($f$)}, and \textbf{codomain}, \emph{cod($f$)}, which are objects in $O(\mathcal{C})$. In general we have
 \[f: \text{\emph{cod($f$)}} \rightarrow \text{\emph{dom($f$)}}\]
 We can also define the \textbf{hom-set} denoted $\mathcal{C}(A,B)$ such that
 \[\mathcal{C}(A,B) := \{f \in A(\mathcal{C}) | f:A \rightarrow B\}. \]
 \item For any triple of objects $A$, $B$, $C$, there exists a \textbf{composition map}
 \[\mathcal{C}_{A, B, C} : \mathcal{C}(A,B) \times \mathcal{C}(B,C) \rightarrow \mathcal{C}(A,C) \]
 This is equivalent to saying for any map, $f$, from $A$ to $B$ and any map, $g$, from $B$ to $C$, there exists a composition map, $g \circ f$ from $A$ to $C$.
 \item For any object, $A$, in $O(\mathcal{C})$ there exists an \textbf{identity} arrow denoted $\text{\emph{id}}_A : A \rightarrow A$.
 \item For any $f: A \rightarrow B$ and two other arrows $g$ and $h$ we have that
 \[\text{\emph{id}}_B \circ f = f = f \circ \text{\emph{id}}_A,\]
 \[h \circ ( g \circ f) = (h \circ g) \circ f. \]
 \end{itemize}
\end{mydef}

\subsection{Diagrams}
The structure of $O(\mathcal{C})$ and $A(\mathcal{C})$ can be easily seen in \textbf{diagrammatic} form. These diagrams show the various associative and identity properties explicity by following the arrows between objects. Say we have a category with the objects $A$, $B$, and $C$, and functions defined as $f : A \rightarrow B$, $g: B \rightarrow C$, and $h: A \rightarrow C$, we can express this category as 

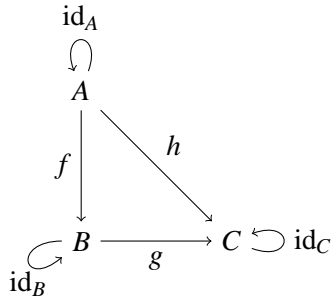
\begin{center}
\begin{tikzpicture}
  \node (A) {$A$}  (A);
  \node (B) [below of=A] {$B$};
  \node (C) [right of=B] {$C$};
  \draw[->] (A) to node {$h$} (C);
  \draw[->] (A) to node[left] {$f$} (B);
  \draw[->] (B) to node [swap] {$g$} (C);
  \draw[->] (B) [out=180, in=220, loop] to node[below] {$\text{id}_B$}(B);
  \draw[->] (A) [out=70, in=110, loop] to node[above] {$\text{id}_A$}(A);
  \draw[->] (C) [out=340, in=20, loop] to node[right] {$\text{id}_C$}(C);
\end{tikzpicture}
\end{center}

In general, we won't include the identity maps in order to have a more presentable diagram

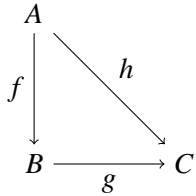
\begin{center}
\begin{tikzpicture}
  \node (A) {$A$} (A);
  \node (B) [below of=A] {$B$};
  \node (C) [right of=B] {$C$};
  \draw[->] (A) to node {$h$} (C);
  \draw[->] (A) to node[left] {$f$} (B);
  \draw[->] (B) to node [swap] {$g$} (C);
\end{tikzpicture}
\end{center}

Nothing prevents us from having the same object or arrow appear multiple times in a diagrammatic description. The property that $\text{id}_B \circ f = f = f \circ \text{id}_A$ can be shown by noting that the following diagram is \textbf{commutative} or path independent

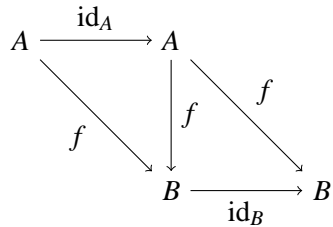
\begin{center}
\begin{tikzpicture}
  \node (A) {$A$} (A);
  \node (A1) [right of=A] {$A$};
  \node (B) [below of=A1] {$B$};
  \node (B1) [right of=B] {$B$};
  \draw[->] (A) to node {$\text{id}_A$} (A1);
  \draw[->] (A) to node[below left]  {$f$} (B);
  \draw[->] (A1) to node {$f$} (B);
  \draw[->] (A1) to node {$f$} (B1);
  \draw[->] (B) to node[below] {$\text{id}_B$} (B1);
\end{tikzpicture}
\end{center}

so that all possible routes from $A$ to $B$ are equivalent (i.e. $f = \text{id}_B \circ f = f \circ \text{id}_A$). We can similarly show associativity with this commutative diagram

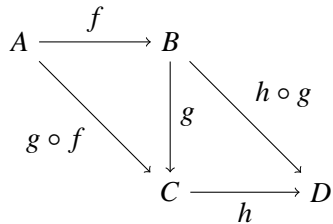
\begin{center}
\begin{tikzpicture}
  \node (A) {$A$} (A);
  \node (B) [right of=A] {$B$};
  \node (C) [below of=B] {$C$};
  \node (D) [right of=C] {$D$};
  \draw[->] (A) to node {$f$} (B);
  \draw[->] (A) to node[below left] {$g \circ f$} (C);
  \draw[->] (B) to node {$g$} (C);
  \draw[->] (B) to node {$h \circ g$} (D);
  \draw[->] (C) to node[below] {$h$} (D);
\end{tikzpicture}
\end{center}

We can also have diagrams based around a square, as in this diagram

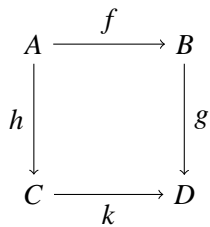
\begin{center}
\begin{tikzpicture}
  \node (A) {$A$} (A);
  \node (B) [right of=A] {$B$};
  \node (C) [below of=B] {$D$};
  \node (D) [left of=C] {$C$};
  \draw[->] (A) to node {$f$} (B);
  \draw[->] (A) to node[left] {$h$} (D);
  \draw[->] (B) to node {$g$} (C);
  \draw[->] (D) to node[below] {$k$} (C);
\end{tikzpicture}
\end{center}

where commuting would mean that $g \circ f = k \circ h$. It is important to note that a diagram doesn't necessarily have to commute unless forced to by the properties of a category (as with the associativity and identity diagrams). Also note that arrows are not in general reversible within a given category.

\begin{mydef}
An \textbf{isomorphism} in a category $\mathcal{C}$ is an arrow $i: A \rightarrow B$ such that there exists an \textbf{inverse} arrow $i^{-1}:B \rightarrow A$ and 
\[i^{-1} \circ i = \text{\emph{id}}_A\; , \;  \;  \;  \;  \; i \circ i^{-1} = \text{\emph{id}}_B. \]
If such an arrow exists, we say that $A$ and $B$ are \textbf{isomorphic}, $A \cong B$.
\end{mydef}

There are also a few specific definitions to categorize possible objects that could exist within a given category. These become important later on as we define special categories which will be applied to foundational logical issues.

\begin{mydef}
An object $I$ in a category $\mathcal{C}$ is \textbf{initial} if, for every object $A \in O(\mathcal{C})$, there exists a unique arrow, $i_A$ from $I$ to $A$.
\end{mydef}

\begin{mydef}
An object $T$ in a category $\mathcal{C}$ is \textbf{terminal} if, for every object $A \in O(\mathcal{C})$, there exists a unique arrow, $t_A$ from $A$ to $T$.
\end{mydef}

\begin{myprop}
If, in category $\mathcal{C}$, there exists two initial elements, $I$ and $I'$, then there exists a unique isomorphism, $\tau : I \rightarrow I'$.
\end{myprop}

\subsection{Products}

The concept of a product is very natural as one looks to see the relations between objects. In traditional set theory, a common construction is a multivariable truth function which takes two objects from a set of all possible statements and maps them to either $0$ or $1$ depending on whether they are mutually exclusive or not. Recall the basic definition in set theory as

\begin{mydef}
The \textbf{Cartesian product} of a set is defined as 
\[A \times B :=\{(a,b)| a \in A \wedge b \in B \} \]
and the set of \textbf{projection functions} as
\[ \chi_1 : (a,b) \rightarrow a , \; \; \; \; \chi_2 : (a,b) \rightarrow b. \]
\end{mydef}

This can easily be generalize in a categorical context. Remember that we still don't assume that there is any structure within the objects, and all the arrows are without any connotation of being mappings of individual elements.

\begin{mydef}
Let $A$,$B$ be objects in $O(\mathcal{C}).$ The object $A \times B$ along with the pair of arrows $\chi_1: A \times B \rightarrow A$, $\chi_2 : A \times B \rightarrow B$ such that for every triple $A$, $B$, $C$, with arrows $f: C \rightarrow A$, $g: C \rightarrow B$ there exists a unique arrow

\[ \langle f,g \rangle : C \rightarrow A \times B \]

where we have the following diagram commuting

\begin{center}
\begin{tikzpicture}
  \node (A) {$A$} (A);
  \node (AB) [right of=A] {$A \times B$};
  \node (B) [right of=AB] {$B$};
  \node (C) [below of=AB] {$C$};
  \draw[->] (AB) to node[above] {$\chi_1$} (A);
  \draw[->] (AB) to node {$\chi_2$} (B);
  \draw[->] (C) to node {$f$} (A);
    \draw[->] (C) to node[below right] {$g$} (B);
    \draw[->] (C) [dashed] to node [above]{$\langle f,g \rangle$} (AB);
\end{tikzpicture}
\end{center}

so that $\chi_1 \circ \langle f,g \rangle = f$ and $\chi_2 \circ \langle f,g \rangle = g$.
\end{mydef}

In order to help solidify these concepts, it is useful to look at a concrete example of a category, the category of sets. An important thing to note when looking at the category of sets is that we are only looking at objects and the arrows between them, and not any sort of structure within the objects. In particular, we do not make mention of elements of the sets except in the context of the functions themselves.\endnote{Lawvere, F. William. ``An elementary theory of the category of sets (long version) with commentary." \emph{Repr. Theory Appl. Categ} 11 (2005): 1-35.}

\begin{myex}
The category of sets, \emph{\textbf{set}}, is a category where all objects $A \in O(\text{\emph{\textbf{set}}})$ are sets,  and all arrows $f$ are defined as total functions between those sets. 
\begin{itemize}
\item Composition of arrows is defined as the normal composition of functions where $g \circ f (a) = g(f(a))$.
\item The identity map $\text{\emph{id}}_A : A \rightarrow A$ is defined as $\text{\emph{id}}_A (a) = a$.
\item Two objects, $A$ and $B$, are isomorphic if all the arrow relations are the same for both sets.
\item The null set, $\varnothing$, is an initial object since one can define a function from $\varnothing$ to an arbitrary set (simply a function without any mappings), but one cannot define a function from a set with elements to the empty set since for a given element there is no object in $\varnothing$ to relate it to.
\item Similarly, singletons (i.e. sets with only one element) are terminal objects as the only function one can construct is that mapping all elements of some set to the single object. All single sets are isomorphic, as one can simply rename the functions by the single object in the set.
\end{itemize}
\end{myex}

\subsection{Functors}

Functors involve two categories, a domain category $A$, and a codomain category $B$, and a mapping which assigns every arrow in $A$ to one in $B$ in a fashion that preserved the commutative properties of a diagram. More formally we have that

\begin{mydef}
A (covariant) \textbf{functor} map $F : \mathcal{C} \rightarrow \mathcal{D}$ is determined by the following

 \begin{itemize}
 \item An object map, assigning an object $FA \in O(\mathcal{D})$ for every $A \in O(\mathcal{D})$.
  \item An arrow map, assigning an arrow $Ff \in A(\mathcal{D})$ for every $f: A \rightarrow B$ in  $A(\mathcal{D})$ such that
  \[ F(g \circ f) = Fg \circ Ff \;\;\;\;\; F\text{\emph{id}}_A = \text{\emph{id}}_{FA} \]
 \end{itemize}
\end{mydef}

There is no reason to believe that a functor need be \emph{bijective} (i.e. that each $FA$ is unique and every element $T \in O(\mathcal{D})$ has a corresponding $B \in O(\mathcal{C})$). In fact, the simplest functor to understand is the \textbf{constant functor}, $F_C$, which maps all objects in a given category to a single object in another category, and all arrows to the identity arrow. Another simple functor is the \textbf{identity functor} which maps all objects to themselves and all arrows to themselves, producing an identical category.

A particularly useful construction in the \textbf{product functor} which is symbolized $- \times Y$. This functor maps all objects $X \in O(\mathcal{C})$ to $X \times Y$, and all arrows $f$ to $\langle f, \text{id}_{Y} \rangle$. Diagrammatically this can be expressed as

\begin{center}
\begin{tikzpicture}
  \node (A) {$A$} (A);
  \node (B) [below of=A] {$B$};
  \node (C) [right of=B] {$C$};
  \draw[->] (A) to node {$h$} (C);
  \draw[->] (A) to node[left] {$f$} (B);
  \draw[->] (B) to node [swap] {$g$} (C);
  \draw[|->,thick] (2.50,-1) -- (3.5,-1);
  \begin{scope}[shift={(5.3,0)}]
  \node (A) {$A \times Y$} (A);
  \node (B) [below of=A] {$B\times Y$};
  \node (C) [right of=B] {$C\times Y$};
  \draw[->] (A) to node {$\langle h, \text{id}_{Y} \rangle$} (C);
  \draw[->] (A) to node[left] {$\langle f, \text{id}_{Y} \rangle$} (B);
  \draw[->] (B) to node [swap] {$\langle g, \text{id}_{Y} \rangle$} (C);
  \end{scope}
\end{tikzpicture}
\end{center}

\subsection{Natural Transformations}
One of the immediate results of the definition of functors is that one can create a category of categories, denoted $\text{\textbf{cat}}$, where all of the objects are categories and all the arrows are functors between them. A level up in abstraction, we can also define a category of functors as objects, and \textbf{natural transformations} as the arrows in the category.

\begin{mydef}
Let $F$ and $G$ be functors $\mathcal{C} \rightarrow \mathcal{D}$. A natural transformation, $t$, is defined as a directional relation from $F$ to $G$, and is indexed by objects in $\mathcal{C}$,

\[ t_A : FA \rightarrow GA \]

or, more schematically,

\[ t: F \rightarrow G. \]
\end{mydef}

If $f$ is an arrow from $A$ to $B$ in category $\mathcal{C}$ we have the following commuting diagram.

\begin{center}
\begin{tikzpicture}
  \node (A) {$FA$} (A);
  \node (B) [right of=A] {$FB$};
  \node (C) [below of=B] {$GB$};
  \node (D) [left of=C] {$GA$};
  \draw[->] (A) to node {$Ff$} (B);
  \draw[->] (A) to node[left] {$t_A$} (D);
  \draw[->] (B) to node {$t_B$} (C);
  \draw[->] (D) to node[below] {$Gf$} (C);
\end{tikzpicture}
\end{center}

The natural transformations between functors between two given categories are the arrows within the \emph{functor category} or \emph{category exponential}. For fixed categories $A$ and $B$, this construction is denoted $B^A$. It can be shown that this object satisfies some nice properties, such as if $t \in A(B^A)$ and $u \in A(C^B)$ then $t \circ u \in A(C^A)$. A useful specific example is the self exponential $A^A$, which essentially describes the structure of the space of endomorphisms (self-mappings).

\subsection{Adjoint}
The property of adjunction plays a crucial role in Lawvere's theories and in the universality of category theory. At its core, adjunction is a duality relationship between functors that satisfies some particular useful and natural properties.

\begin{mydef}
An \textbf{adjoint situation} involves categories ($A$ and $B$), two functors ($F: A \rightarrow B$ and $G: B \rightarrow A$), and two natural transformations ($t : A \rightarrow FG \in A(A^A)$ and $u: GF \rightarrow B \in A(B^B)$) such that the following two relations hold,

\[ tF\circ Fu = F, \; \; \; \; \; \; Gt \circ uG = G. \]

$F$ and $G$ are said to \textbf{adjoint} to one another. Adjoint situations are sometimes listed by their functors and natural transformations as $\langle G, F, t, u \rangle$.
\end{mydef}

It can be shown\endnote{Lawvere, F. William. ``Adjointness in foundations." \emph{Dialectica} 23.3-4 (1969): 281-296.} that adjoint situations have some important properties.

\begin{itemize}
\item Adjoints are unique, in that given any two adjoint situations with a given functor $F$, the adjoint functors $G$ and $G'$ will be isomorphic.
\item Given two fixed categories $A$ and $B$ in two adjoint situations, any natural transformation $v : F \rightarrow F'$ induces a unique natural transformation $w: G \rightarrow G'$.
\item Similarly adjoint situations can be composed and remain adjoint. If one has two adjoint situations , $\langle G, F, t, u \rangle$ between $A$ and $B$ and $\langle G', F', t', u' \rangle$ between $B$ and $C$, one can construct 
\[ F' \circ F : A \leftarrow C \]
which left adjoint to 
\[ G \circ G' : A \rightarrow C. \]
\end{itemize}

%
%

\subsection{Cartesian Closed Categories}

Cartesian closed categories have additional structure imposed on the relationship of arrows within the category and the existence of exponential objects. They have found widespread applications in logic, not only in the context of categorical constructions of Tarski's theorem, but also in formulations of lambda-calculus\endnote{Huet, Gerard. ``Cartesian closed categories and lambda-calculus." \emph{Combinators and Functional Programming Languages.} Springer Berlin Heidelberg, 1986. 123-135.} and theoretical computer science more generally.\footnote{Often-times Lawvere's work refers to a \emph{topos}, which is a Cartesian closed category with additional substructure conditions. We will limit our discussion to just Cartesian closed categories.}

\begin{mydef}
A \textbf{closed Cartesian category} is category $\mathcal{C}$ which satisfies the following properties:\footnote{Lawvere defined closed Cartesian categories in terms of the adjoint relationships within them which result in the listed properties.}
\begin{itemize}
\item Any two objects, $A,B \in O(\mathcal{C})$, have a product $A \times B \in O(\mathcal{C})$.
\item  Any two objects, $A,B \in O(\mathcal{C})$, have a exponential $B^A \in O(\mathcal{C})$
\item There exists a terminal object, denoted $\mathbb{1}$.
\end{itemize}
\end{mydef}

The category \textbf{set} is cartesian closed since the product set is just the cartesian product of two sets (a set) and the exponential, $B^A$ is the set of all functions from $A$ to $B$, which is also a set.\endnote{Lambek, Joachim. ``Deductive systems and categories III. Cartesian closed categories, intuitionist propositional calculus, and combinatory logic." \emph{Toposes, algebraic geometry and logic.} Springer Berlin Heidelberg, 1972. 57-82.}

\section{Category Theory vs. Set Theory}

\indent 
In considering the foundations of mathematics, the most common approach has been couched  in terms of a set-theoretic axiomatization of the ``member relationship".  It is a remarkable empirical fact that all of modern mathematics can be placed within some formulation of set theory (although different concepts may need different set theoretic formulations). To be more precise, all mathematical objects can be considered sets through a hierarchical structure, and all the crucial properties of such sets can be derived from the axioms of a set theory.\endnote{Maddy, Penelope. \emph{Defending the axioms: On the philosophical foundations of set theory.} Oxford: Oxford University Press, 2011.} It seems perfectly plausible that one could look just as easily at axiomatization of the notion of a ``function"  as opposed to a set.\footnote{An axiomatized fuctional formulation is important in $\lambda$-calculus, for example.} The categorical notion of arrow is merely an extension of the notion of a function to simply be a directional relationship, with any associated structure additional for that particular case (such as for \textbf{Set}).

\indent A natural criticism of category theory is that it doesn't provide a true divergence from the line of set theoretical constructions. Tarski himself asked at a conference, ``But what is a category if not a set of objects together with a set of morphisms?" Lawvere reply stressed that foundational set theory is based on the binary relationship of membership, while category theory is based on a ternary relationship of composition (i.e. commutativity of diagrams). Category theory in this way focuses on structure, while set theory focuses on identity. A single object of category theory is effectively meaningless, while a single set in set theory has its own inherit meaning through concepts of membership, cardinality, etc. Relations between objects are how structure emerges in category theory; specifically one identifies a codomain of a given arrow, thereby allowing to explicitly differentiate between an inclusion map and an identity map (i.e. set theory doesn't distinguish fully between $f: \mathbb{R} \rightarrow \mathbb{R}$ and $f: \mathbb{R} \rightarrow \mathbb{N}$).

\subsection{Set Theoretic Approach to Category Theory}

\indent  The question still remains what is (or if there exists) the appropriate set theoretical foundation category theory. There are three commonly suggested answers to this question; Grothendieck Universes, reflection principles, or none whatsoever. Grothendieck universes are perhaps the most natural of these to study category theory as they provide an exact list of properties that keep the sets (and therefore categories) ``small", in such a way to avoid the failings of naive set theory. 

\begin{mydef}
A Grothendieck universe is a set $U$ with the following properties
\begin{itemize}
\item If $x,y \in U$ then $ \{x,y\} \in U$.
\item if $x \in U$ and $y \in x$ then $y \in U$.
\item If $x \in U$ then the power set of $x$ is also in $U$.
\item if the exists a family of sets, $\{x_a\}_{a \in I}$, within $U$, and if $I \in U$ then the union of all $x_a$ is also an element of $U$.
\end{itemize}
\end{mydef}

This approach to construction has some immediate drawbacks. It implies the existence of an inaccessible cardinal, $k$, and all mathematics are formed before stage $k$ (i.e. in the process of taking unions of sets). One cannot construct a category of categories, which is a fairly natural process in category theory and other ``large" sets. We also have to extend beyond standard ZF(C) constructions to establish this inaccessible cardinal.\footnote{The restricted environment that Grothendieck Universes are constructed is often known as $\text{ZF+}$.}

Reflection principle has a similar hierarchical structure but attempts to avoid the issues of inaccessible cardinals by saying that the universe of small sets is simply an elementary substructure of the universe of all sets. More precisely, we have that each first-order statement within the theory of small sets has the same meaning when the variables range over all sets. With this in mind, we only need to prove categorical theorems within the small-set structure and then ``reflect" that result to all sets. We also need not extend beyond ZF(C) set theory to construct the basic relations to be reflected.\endnote{Reinhardt, W. N. (1974), ``Remarks on reflection principles, large cardinals, and elementary embeddings.", \emph{Axiomatic set theory}, Proc. Sympos. Pure Math., XIII, Part II, Providence, R. I.: Amer. Math. Soc., pp. 189Ð205,}

The effort to develop category theory entirely without the structure of sets is very much a goal of those who view category theory as a ``better" foundation to mathematics or, at least, a useful tool distinct from set theory. The primary promotor of this view is Saunders MacLane, the cofounder of category theory, who formulated the concept of a \emph{metacategory} entirely within a first order system without any concept of ``set" built in.\endnote{Mac Lane, Saunders. \emph{Categories for the working mathematician.} Vol. 5. Springer verlag, 1998.} He avoids the use of the concept of a set of objects and the terminology associated with it, instead focusing purely on relational aspects between objects. According to him, there are two fundemental concepts of category theory that standard ZF set theory cannot describe; (i) One cannot form a category of all structures of a given kind (i.e. category of sets, category of categories) due to self referential paradoxes, and (ii) one cannot form the exponential category $B^A$ of any two given categories, but doesn't himself offer a consistent way to construct these objects without resorting to some set-theoretic notion of size.

Utilizing the concept of a fibration (a fundamentally leveled structure), others have attempted to find a new way \emph{within} category theory to define size to avoid self-referential paradoxes.\endnote{Benabou, Jean. ``Fibered categories and the foundations of naive category theory." \emph{The Journal of Symbolic Logic} 50.1 (1985): 10-37.} It is not surprising that the reliance of most set-theoretical foundational constructions on a constructed hierarchy of cumulative types might also be found in a the construction of foundational theories utilizing categories.

\subsection{Categorical Approach to Set Theory}

There have been various attempts to construct axiomatized set theory through a more structured version of a cartesian Closed category with a ``truth object" incorporated in it. Through various possible arrows to the truth object, it has been suggested that one can construct any possible axioms for a set theory implicitly.\endnote{Awodey, Steve. ``From Sets to Types, to Categories, to Sets." \emph{Foundational Theories of Classical and Constructive Mathematics. }Springer Netherlands, 2011. 113-125.}

It has also been shown that very fundamental theorems generally stated through set theory, such as Tarski's Undefinability Theorem and G\"{o}del's incompleteness theorems, as well as statement of paradoxes, can be cast completely within the framework of category theory as a natural consequence of the Fixed Point Theorem which limits possible arrows between objects in a closed Cartesian category.\endnote{Lawvere, F. William. ``Diagonal arguments and cartesian closed categories." \emph{Category theory, homology theory and their applications II. }Springer Berlin Heidelberg, 1969. 134-145.} This common structure between seemingly disjoint paradoxes and theorems suggests that there could be deep insight found in category theory when viewed as a tool in mathematical foundations and a simple extension of category theory to set theory.\endnote{Pavlovic, Dusko. ``On the structure of paradoxes." \emph{Archive for Mathematical Logic} 31.6 (1992): 397-406.} 

More contentiously, one can simply view a category as a form of a generalized set. This is perhaps a natural assumption to make, as \textbf{Set} is suggested to simply be one of many possible categories. If there is success in constructing category theory of such large categories without resorting to set-theoretical notions then set theory has no particular priority within the framework of mathematics beyond possible psychological implications concerning the primacy of the elemental notions of set theory. However, so far, many logicians seem unconvinced of the ability of category theory to establish the notion of size within itself. The table below presents the objects in set theory as possible restricted objects in category theory.

\begin{center}
\begin{tabular}{ |l|l| }
  \hline
  \multicolumn{2}{|c|}{Category Theory vs. Set Theory} \\
  \hline
  Object & Set \\
  Morphism & Function \\
  Monomorphism & One-to-one function\\
  Epimorphism & Surjection \\
  Isomorphism & Bijection\\
  Product & Cartesian product \\
  Coproduct & Disjoint union\\
  \hline
\end{tabular}
\end{center}
%
%

\section{Universal Constructions}

\indent Category theory provides an extremely useful framework to relate different mathematical concepts. Although originally developed in the context of algebraic topology, various objects constructed in category theoretic terminology (such as cartesian products, power sets, universal covering spaces, Stone-Cech compactification, etc.) can easily be applied to other fields of mathematics. The concepts of universal constructions and the relationships between them within seemingly different fields, has lead to many innovative proofs for otherwise complex theorems in algebraic topology and algebraic geometry. In this way, universal constructions allow one to focus on underlying structure as opposed to complicated technical aspects of a given field.

The most commonly exploited aspect of category theory is the general usefulness of an adjoint operators. The uniqueness results mentioned in the previous section have extremely important implications in constructing given mappings within various settings. For example, one can see through universal structures that free groups, universal enveloping algebras and Stone-Cech compactifications are all simply the left adjoint of a functor.\endnote{Blass, Andreas. ``The interaction between category theory and set theory." \emph{Preprint, University of Michigan} (1984).}

Category Theory's role as a universal construction is one of the foundational points of Kreisel's objection to viewing category theory as a ``foundational tool." He distinguishes between the concept of a foundational structure and an organizational structure. A foundational structure (in his view, axiomatized set theory) has to be ``as specific as possible" in order to ``make strong assertions." Foundations should be expected to create structures that do not occur in practice, like type theory of sets, while practical application in the proofs of everyday mathematicians do not require distinction in this regards. Universal constructions are an example of a less-specific construction which doesn't play a role in foundational specified constructions. In this way they are useful for \emph{practice} of mathematics and increasing the ``efficiency of proofs", not in the ``analysis of proofs." Throughout his discussion, he sees foundational (set theory) being fundamentally distinct from organizational (category theory).\endnote{Kreisel, Georg. ``Observations on popular discussions of foundations." \emph{Axiomatic set theory}. No. 13 part 1. 1971.}

However, I view this distinction, particularly in regards to category theory and set theory's role in the foundational mathematics as misguided. Throughout mathematics one has essentially two dual aspects; the Formal and the Notional. In basic algebra, we have a distinction between the manipulation of a polynomial function and the visualization of the resulting graphical curve. We have the set of theorems that govern a class of mathematical objects and then we have the analysis of the specific objects in question. In these cases we see that the Notional is the subject matter of the Formal. In the same way, when one considers foundational mathematics there could exist a similar, unexplored, relationship between sets and categories. Although each can provide insight to one another, set theory has a formal component as explicitly defining its constructions, while category theory is more notional in its search for generic structures. To say that one is more foundational than the other ignores the mutual dependancies which are being developed.


\section{Conclusion}

\epigraph{Beauty is the first test: there is no permanent place in the world for ugly mathematics.}%
         {Godfrey H. Hardy (1877-1947) \\ A Mathematician's Apology}

When considering the foundations of mathematics, it is clearly difficult to fully reconcile category theory and set theory. While much work has been done to axiomize category theory to such an extent that it can stand on its own, in practice one still needs to introduce the concepts of size from a given set theoretic standpoint. On the other hand, set theory seems incomplete in fully grasping nuances of category theory and it is unable to make relations between objects in mathematics as explicitly clear as category theory. Although there are no notable ``proofs" that category theory establishes that a given axiomatization of set theory cannot, category theory has shown itself to be able to expose fundamental properties more so than any other mathematical sub-field. Until category theory is able to fully encompass set theory, it would seem remise to ignore either when considering fundamental questions at the core of mathematical inquiry.

\theendnotes
\end{document}